\documentclass[11pt]{amsart}
\usepackage{amssymb}
\usepackage[
            pdftex,
           colorlinks=true,
            urlcolor=blue,       
            filecolor=blue,     
            linkcolor=blue,       
            citecolor=blue,         
            pdftitle= {Effective reconstruction of generic genus 5 curves from
theta hyperplanes},
            pdfauthor={David Lehavi},
            pdfsubject={},
            pdfkeywords={}
            pagebackref,
            pdfpagemode=UseThumbs,
            bookmarksopen=true]
            {hyperref}

\newcommand{\BP}{{\mathbb{P}}}
\newcommand{\BA}{{\mathbb{A}}}
\newcommand{\BC}{{\mathbb{C}}}
\newcommand{\BZ}{{\mathbb{Z}}}
\newcommand{\BF}{{\mathbb{F}}}
\newcommand{\ii}{{\mathfrak{i}}}
\newcommand{\bW}{{{W}}}

\newcommand{\gd}{\delta}
\newcommand{\gb}{\beta}
\newcommand{\gc}{\gamma}
\newcommand{\gC}{\Gamma}
\newcommand{\gS}{\Sigma}

\newcommand{\gep}{\epsilon}
\newcommand{\ga}{\alpha}
\newcommand{\gth}{\theta}
\newcommand{\gT}{\Theta}

\newcommand{\cO}{{\mathcal{O}}}
\newcommand{\cV}{{\mathcal{V}}}
\newcommand{\cL}{{\mathcal{L}}}

\newcommand{\ol}[1]{\overline{#1}}

\newcommand{\Hess}{\mathrm{Hess}}

\newcommand{\spn}{\mathrm{span}}

\newcommand{\Pic}{\mathrm{Pic}}

\newcommand{\Sym}{\mathrm{Sym}}
\newcommand{\Div}{\mathrm{Div}}
\newcommand{\sm}{\smallsetminus}
\theoremstyle{plain}
\newtheorem{lma}{Lemma}[section]
\newtheorem{thm}[lma]{Theorem}

\newtheorem{prp}[lma]{Proposition}
\newtheorem{cor}[lma]{Corollary}
\theoremstyle{definition}
\newtheorem{prd}[lma]{Proposition-Definition}
\newtheorem{dfn}[lma]{Definition}
\newtheorem{rmr}[lma]{Remark}
\newtheorem{ntt}[lma]{Notation}
\newtheorem{dsc}[lma]{}

\begin{document}
\title[Reconstruction of genus 5 curves from theta hyperplanes]{Effective reconstruction of generic genus 5 curves from their theta hyperplanes}
\begin{abstract}
We effectively reconstruct the set of enveloping quadrics of a generic curve
$C$ of genus $5$ from its theta hyperplanes; for a generic genus $5$ curve $C$
this data suffices to effectively reconstruct $C$.
As a consequence we get a complete description of the Schottky locus in genus
$5$ in terms of theta hyperplanes.
The computational part of the proof is a certified numerical argument.
\end{abstract}
\author{David Lehavi}
\email{dlehavi@gmail.com}
\date{\today}
\maketitle
\setcounter{tocdepth}{1}
\tableofcontents
%
\section{Introduction}
%
The quest for methods of reconstructing a curve from its theta hyperplanes
goes back to
the 19th and early 20th century geometers Aronhold and Coble (see
\cite{A}, \cite{Co} chapter IV, and \cite{D} sections 6.1.2 and 6.2.2):
in the non
hyperelliptic genus 3 case, theta hyperplanes are simply
bitangents, and both Aronhold and Coble provided formulas for reconstructing
curves from certain ordered subsets of the $28$ bitangents of the curve.

Recent years witnessed some revived interest in generalizations of this
problem from several directions: first relaxing the need for {\em ordered}
theta hyperplanes (see \cite{CS1}, \cite{L1}) (where ``order'' means
known action of the level-2 cover deck transformations group $SP_{2\mathrm{genus}}(2)$),
and then generalizations to
higher genus curves (see \cite{CS2}) and Abelian varieties (see \cite{GS-M1},
\cite{GS-M2}). Using the fact that \cite{CS2} relaxed the need for ordered theta
characteristic, we gave in \cite{L2} two results for
generic genus $4$ curves: An explicit construction of the enveloping quadric,
and an explicit construction of the nodes of Wirtinger cubics, which \cite{CKRN}
used to construct the variety of enveloping cubics of $C$, thus completing
the reconstruction of $C$.

The object of this paper is to extend --- in Theorem~\ref{T:main} --- the
first of the two results from \cite{L2} to genus $5$.
Throughout this paper we consider a generic complex curve $C$ of
genus $5$; since $C$ is generic we assume that all its odd theta
characteristics
are $1$ dimensional - i.e.\ if $\gth$ is an odd theta characteristic of $C$
then $\dim H^0(C,\gth)=1$.
Hence, for each
odd theta characteristic $\gth$ there exists a unique hyperplane $l_\gth$
--- called a {\em theta hyperplane} --- in the
dual canonical system of $C$ such that when $C$ is identified with its
canonical image, the points in the intersection product
$C\cdot l_\gth$ are all double, and the points in $\frac{1}{2}C\cdot l_\gth$
sum up to $\gth$. This hyperplane is the projectivization of the plane
$T_\gth\gT_C\subset T_\gth J C$ under the identification of $T_\gth J C= T_0 J C= H^0(K_C)^*$, and where $J C,\gT_C$ are the Jacobian variety and theta divisor of
the curve $C$.
Recall (see e.g.\ \cite{D} 5.4.2), that if $\ga$ is a non-trivial
$2$ torsion point on the Jacobian $J C$, then the
{\em Steiner system} $\gS_{C,\ga}$ of the pair $(C,\ga)$ is defined to be the
set
\[
\{\gth: 2\gth = K_C\text{ and } \dim H^0(C,\gth)=\dim H^0(C,\gth+\ga)\equiv 1 \mod 2\}.
\]
The number of theta characteristics in a Steiner
system of a genus $g$ curve is $2^{g-1}\cdot(2^{g-1}-1)$; i.e.\ in our case
a Steiner system is comprised of $16\cdot15=240$ odd theta characteristics,
out of the total of $2^{5-1}(2^5-1)=496$
odd theta characteristics of the curve.
For each pair $\gth,\gth+\ga$, and corresponding theta hyperplanes
$l_\gth, l_{\gth+\ga}$, we let $q_{\{\gth,\gth+\ga\}}\in |\cO_{|K_C|^*}(2)|$ be the
image of $\{l_\gth, l_{\gth+\ga}\}\in S^2|K_C|$ under the map:
$S^2|K_C|\to |\cO_{|K_C|^*}(2)|$. We can now state the following:
\begin{thm}\label{T:main}
  Let $C$ be a generic genus $5$ curve, with $q_{\gth, \gth+\ga}$ as above,
  then the set of quadric enveloping $C$ is given by:
  \[
  I_2(C)=\bigcap_{\ga\in J C[2]\sm\{0\}}\spn(\{q_{\{\gth,\gth+\ga\}}\}_{\gth\in\gS_{C,\ga}})
  \]
\end{thm}
We note two immediate corollaries of this theorem:
\begin{cor}
  For a generic genus $5$ curve $C$, Theorem~\ref{T:main} effectively
  determines the curve.
\end{cor}
\begin{proof}
  By the Enriques-Babbage Theorem (see e.g.\
  \cite{ACGH} chapter VI \S3), a generic canonical curve is either
  trigonal, isomorphic to a plane quintic, or is
  cut out by quadrics, in which case $C=\cap_{[q]\in I_2(C)}Z(q)$.
\end{proof}
Since theta hyperplanes are defined for any principally polarized
Abelian variety $(A,\gT)$ via the Gauss map, we may define
$\mathbb{I}_2(A, \gT)$ as the right hand side
of the equation in Theorem~\ref{T:main}. Using this definition we get:
\begin{cor}
  The Schottky locus in genus $5$ is comprised of (the closure of) pairs
  $(A,\gT_A)$ such that the intersection $\cap_{[q]\in \mathbb{I}_2}Z(q)$ is a
  curve $C$ of arithmetic genus $5$, and
  such that the moduli points $[(J C, \gT_C)],[(A,\gT)]$ in the moduli of
  principally polarized Abelian varieties are identical.
\end{cor}
\begin{proof}
  This follows from Theorem~\ref{T:main} and the previous corollary.
\end{proof}
The paper is organized as follows: In Section~\ref{S:special} we show how to
reduce the proof of the theorem to slightly stronger claims on a specific curve:
as in \cite{L2} most of the proof is done by specializing to one curve; in
the case at hand the Wiman curve $\bW=W^{160}$.

Sections~\ref{S:humbert}, \ref{S:steiner} and \ref{S:rep} form the bulk of the
algebro-geometric
part of the paper. In Section~\ref{S:humbert} we recall the definition of
Humbert curves,
and concretely describe a set of $160$ theta characteristics on them. In
the
remainder of the paper we carefully analyze this set of odd theta
characteristics. In
Section~\ref{S:steiner} we give a partial description of the intersections of
pairs of points in the set of
$160$ theta characteristics we found with the Steiner systems of Humbert
curves:
we show that these pairs intersect exactly $510$ out of the
$1023$ Steiner
systems, and compute the sizes of the intersections.
In Section~\ref{S:rep} we introduce the Wiman curve as a special
Humbert curve with a large automorphism group. We compute
explicit coordinates for the theta characteristics we found in Section\
\ref{S:humbert}, and decompose the second symmetric tensor of the canonical
system of $\bW$ to irreducible representation of its automorphism group.
The most immediate gain of the later is a natural decomposition of the space of
quadrics in the canonical system
to a direct sum of $I_2(\bW)$ and $H^0(2K_\bW)$.

In Section~\ref{S:approx} we show how to produce a certified numerical proof
that a given sum of four theta characteristics on $\bW$ is not 2-canonical. In
a nutshell, we show how to certify that all the
non-trivial quadrics in $H^0(2K_\bW)$ --- considered as the orthogonal
complement of $I_2(\bW))$ --- are bounded away from zero on the sum of the
theta characteristics at hand. We note that while certified numerical arguments
are not new to algebraic geometry (see e. g. \cite{BST}, \cite{HS}, and the
\verb|bertini| package),
they usually deal with a harder case of assuring the existence of a true zero in
the neighborhood of an approximately computed one. The approach we take is
closer in spirit to e.g.\ \cite{OS}; we simply bound values {\em away} from zero.

Finally, in Section~\ref{S:prog} we connect the results from Section\
\ref{S:approx} with actual numerical error bounds, tie it with the computation
of Section~\ref{S:steiner} to completely describe the intersections of
pairs of theta characteristics found in Section~\ref{S:humbert} with all Steiner
systems; then use the decomposition
into irreducible representations from Section~\ref{S:rep} to produce a
certificate for the numerical computation of the
dimension of the intersection in Theorem~\ref{T:main}. We
also provide a witness (in the proof theoretic sense), which can be verified
by a relatively small number of ``traditional'' matrix rank computation
over a bi-quadratic extension of the rational numbers. It is important to note
that the computer program accompanying this paper, and
described in Section~\ref{S:prog} is an integral part of this paper. This
program can be found at {\small \verb|https://arxiv.org/src/1908.02355v2/anc/W160_steiner_systems_and_IC2.cc|}.

Sadly, unlike the theorem itself, which I conjecture
generalizes to higher genera, it is hard to see how the proof technique from
this paper and from
\cite{L2} will generalize to genera greater than $14$, where one cannot
work directly with a ``general curve''. In the ``next'' cases, genera 6,7,
the obvious candidates to work with are curves with large automorphism groups:
e.g.\ some curve along the Wiman-Edge pencil (see \cite{DFL}) in genus 6,
and the Fricke-Macbeath curve in genus 7 (see \cite{M}). We
note the genus 6 case is especially interesting for two reasons: it is the
first genus
where the general curve is not a complete intersection, and it is the first
case where the Schottky problem is not fully resolved.
\subsection*{Acknowledgments}
I would like to thank T. Celik and A. Kulkarni for many useful
remarks they gave on an early draft of this paper.
%
\section{Proving the general result by specialization}\label{S:special}
%
\begin{ntt}
  Given a curve $C$ and a divisor $D\in\Div C$, denote by $\phi_D$
  the map from $C$ to the dual linear system $|D|^*$;
  we identify $C$ with it's canonical image $\phi_{K_C}(C)$.
  Denote by $\ga$ some point in $J C[2]\sm\{0\}$.
  Denote the null set of a homogeneous polynomial $f$ by $Z(f)$.
  Finally denote the projectivization of a space $V$ by $\BP V$.
\end{ntt}
\begin{dfn}[Steiner spaces]\label{D:steiner_sys}
  We set
  \[
  \begin{aligned}
    p_D:H^0(\cO_{|D|^*}(2))\sm\{q:q|_{\phi_D(C)}=0\}&\to H^0(2D)\\
      q&\mapsto Z(q)\cdot \phi_D(C).
  \end{aligned}
  \]
  Note that there are natural isomorphisms
  $H^0(\cO_{|D|^*}(2))\cong\Sym^2 H^0(D)$.
  With $C,\ga$ as above, denote by $V_{C,\ga}$ the closure of the pre-image
  \[
  p_{K_C}^{-1}p_{K_C+\ga} H^0(\cO_{|K_C+\ga|^*}(2))\subset H^0(\cO_{|K_C|^*}(2)).
  \]
  The authors of \cite{CKRN} suggested to name $V_{C,\ga}$
  {\em Steiner spaces}.
\end{dfn}
We want to reduce Theorem~\ref{T:main}, which is stated for a generic curve, to
a somewhat stronger claim about a single specific curve. To this end we use
a standard degeneration argument:
\begin{lma}[see e.g.\ Corollary\ 7 from \cite{L2}]\label{L:C7L2}
  Let $\cV/X$ be a vector bundle over a base $X$, and let
$\cV_1,\dots \cV_n$ be sub-bundles of $\cV$.
Then the function $\dim \langle \cV_1|_x,\ldots, \cV_n|_x\rangle$
is lower semi-continuous on $X$,
and the function
\[
\dim (\cap_{i=1}^n\cV_i|_x)
\]
is upper semi-continuous on $X$.
\end{lma}
Armed with the lemma we now prove the following reduction:
\begin{prp}\label{P:specialize}
  To prove Theorem~\ref{T:main} it suffices to exhibit one non-hyperelliptic
  curve $C_0$ and a set $A\subset J C_0[2]\sm\{0\}$ so that
  \[
  \dim\spn(\{q_{\{\gth,\gth+\ga\}}\}_{\gth\in\gS_{C_0,\ga}})=13,\quad\text{ for all }
  \ga\in A,
  \]
  and
  \[
  \dim\cap_{\ga\in A}\spn(\{q_{\{\gth,\gth+\ga\}}\}_{\gth\in\gS_{C_0,\ga}})=3
  \]
\end{prp}
\begin{proof}
  It is clear from Definition \ref{D:steiner_sys} that
  \[
  V_{C,\ga}\supset\spn(\{q_{\{\gth,\gth+\ga\}}\}_{\gth\in\gS_{C,\ga}}), I_2(C),
  \]
  and that
  \[
  \dim V_{C,\ga}\leq \dim\Sym^2 H^0(K_C+\ga) + \dim I_2(C).
  \]
  Thus, to prove that $V_{C,\ga}$ and
  $\spn(\{q_{\{\gth,\gth+\ga\}}\}_{\gth\in\gS_{C,\ga}})$
  are equal it suffices to show that for the curve $C_0$:
  \[
  \begin{aligned}
    \dim\ &\spn(\{q_{\{\gth,\gth+\ga\}}\}_{\gth\in\gS_{C_0,\ga}})=
    \dim\Sym^2 H^0(K_{C_0}+\ga) + \dim I_s(C_0)
    \\ =&{\dim H^0(K_{C_0}+\ga) + 1\choose 2} + (\dim \Sym^2 H^0(K_{C_0}) -
    \dim H^0(2K_{C_0}))
    \\ =& {4 + 1\choose 2} + {5 + 1\choose 2}- (4(5-1)-(5-1))=13.
  \end{aligned}
  \]
  The claim now follows from Lemma~\ref{L:C7L2}
\end{proof}
%
\section{Odd theta characteristics of Humbert curves}\label{S:humbert}
%
\begin{prd}[Humbert's curves, originally defined in \cite{H}]
  Assume $C$ is a non hyperelliptic genus $5$ curve, then $C$ sits
  on a smooth quadric if and only if is it is non trigonal (see \cite{ACGH},
   chapter VI exercises F1, F2).
  Henceforth we further assume that $C$ is also non-trigonal.
  Denote by $\gC\subset\BP I_2(C)$ the locus of quadrics of rank $\leq4$, and
  by $\gC'$ the locus of quadric of rank $\leq 3$. Denote
  by $W^r_d(C)\subset\Pic^d(C)$ the subscheme of degree $d$ divisors $D$ so that
  $\dim |D|\geq r$.
  Then by \cite{ACGH} chapter VI exercise batch F, the following properties
  hold (for proofs see the specific exercises in each parenthesis):
  \begin{itemize}
  \item $\gC$ is a plane quintic with no multiple components and at most
    ordinary nodes, which are exactly $\gC'$ (see exe. 4).
  \item The map $\chi:W^1_4(C)\to\gC$ defined by
    $D\mapsto\BP(\text{Tangent cone}_D W^0_4(C))$
    is a double cover, branched over $\gC'$. Specifically,
    $W^1_4(C)_{\text{sing}}=\{D\in W^1_4(C):2D=K_C\}$ (see exe. 5, 6, 7).
  \item $W^1_4(C)$ is irreducible if and only if $\gC$ is (see exe. 10).
  \item If $C\to E$ is a bielliptic double cover, then there is a
    component $\gS$ of $W^1_4(C)$ corresponding to the $g^1_2$s on $E$.
    These components project to line components on $\gC$, and vice-versa.
    (see exe. 11, 12).
  \item $C$ has five bi-elliptic pencils if and only if it has 10
    semi-canonical pencils, if and only if, after a suitable change of
    coordinates --- called {\em diagonalization} --- $\BP I_2(C)$ is
    spanned by three quadrics of the form $Q_\ga=\sum_{i=0}^4\lambda_{\ga,i}x_i^2$.
    (see exe. 13, 14). These are called {\em Humbert curves}.
  \end{itemize}
  Collecting the properties above we see that on a Humbert curve
  the variety $W^1_4(C)$ consists of five subvarieties $\Pic^2(E_i)$ ---
  where each
  $E_i$ is an elliptic double quotient of $C$ --- which intersect at
  $\binom{5}{2}=10$ points.
  Each of these five $\Pic^2(E_i)$ is a double cover of one line in $\gC$ under
  the $\chi$.

  For an alternative presentation to the one above, see either \cite{H},
  \cite{E} or \cite{V}.
\end{prd}
\begin{prd}\label{p:elliptic_cone}
  Let $C$ be a Humbert curve, and assume diagonalized coordinates
  $x_0,\ldots, x_4$ on the
  canonical system. By elimination on the monomials $x_0^2,\ldots, x_4^2$
  there are exactly $\binom{5}{3}= 10$
  quadrics (up to constants) in $I_2(C)$ supported on only three coordinates.
  Picking any four coordinates, four of these
  10 quadrics are supported on a subset of the four coordinates. By definition,
  the points corresponding to these four quadrics in $\BP I_2(C)$ sit on a line,
  and
  the intersection of their null sets is a cone over an elliptic curve (the
  intersection of the quadrics in 3-dimensional projective space) in the
  canonical system of $C$. Thus we have identified the elliptic double
  quotients $E_i$; we fix the indexing so that $i$ is the
  vanishing coordinate of the cone over $E_i$.
  As for the realizing the semi-canonical pencils, each of the
  10 quadrics supported on 3 coordinates is a cone over some conic in the
  $\BP^2$ generated by these coordinates.

  Taking indices modulo $5$, denote the semi-canonical pencil of
  the form
  $ax_{i-1}^2+bx_i^2 +cx_{i+1}^2$ by $Q_i$, and the semi-canonical pencil of the
  form $ax_{i-2}^2+bx_i^2 +cx_{i+2}^2$ by $Q_i'$. Denote by $\pi_i$ the five double
  elliptic covers $C\to E_i$. We will identify $\Pic^2E_i$
  with its image in $W^1_4(C)$.
  Denote the points representing $Q_i,Q_i'$ on $\gC'\subset\gC\subset I_2(C)$
  by $[Q_i],[Q_i']$.
  Using these notations we have
  \[
  \chi(\Pic^2E_i)\cap\gC'=\{[Q'_{i+1}], [Q'_{i-1}], [Q_{i-2}], [Q_{i+2}]\}.
  \]
  Denote
  \[
  \gb_i:=\chi^{-1}[Q_{i+2}]-\chi^{-1}[Q_{i-2}],\quad
  \ga_i:=\chi^{-1}[Q_{i+2}]-\chi^{-1}[Q'_{i-1}].
  \]
\end{prd}
Henceforth, we will always assume our Humbert curves are diagonalized, with
$E_i, Q_i, Q_i'$ as in Proposition-Definition \ref{p:elliptic_cone}.
\begin{prp}
  Consider $\ga_i, \gb_i$ as degree 0 divisors on $\Pic^2E_i\subset W^1_4(C)$,
  then
  any curve isomorphism $\Pic^2 E_i\to E_i$, induces the same natural
  isomorphism $\Pic^0(\Pic^2 E)=\Pic^0 E$. Moreover, this isomorphism takes
  $\ga_i, \gb_i$ to distinct non-trivial points in $\Pic E_i[2]$
\end{prp}
\begin{proof}
  The first claim follows since for any elliptic curve we have
  $\Pic^0 E\cong E$ as curves. As for the second claim,
  they are differences between one ramification
  point under $\chi$, and two other ramification points.
\end{proof}
\begin{ntt}
  Henceforth we identify $\ga_i, \gb_i$ with their respective images in
  $\Pic^0(E_i)$.
\end{ntt}
\begin{prp}[The kernel of the map $\prod E_i \to J C$]
  The map $\prod J E_i\to J C$ is a degree $2^5$ isogeny of principally
  polarized Abelian varieties,
  whose kernel is generated by $\{\gb_i-\ga_{i+1}-\ga_{i-1}\}_{i=1}^5$, where
  indices are taken mod $5$.
\end{prp}
\begin{proof}
  The inverse image of the point $(e_0, e_1, e_2, e_3, e_4)\in\prod\Pic^1E_i$
  under the map $\prod_i\pi_i:\Pic C\to \prod_i\Pic E_i$ is $\{(c_0,\ldots,c_4)|c_i\in\pi^{-1}e_i\}$.
  Hence the map
  $J C\to\prod J E_i$ is of degree $2^5$, hence so is the dual map
  $\prod J E_i \to J C$. We now consider, for a given $i$, the points
  $\chi^{-1}[Q_{i+2}], \chi^{-1}[Q_{i-2}], \chi^{-1}[Q'_i]$.
  Computing on $E_i$, we see that the differences between the first two points
  is the image of $\gb_i$ under
  the map $\prod \Pic E_i\to \Pic C$. Computing on $E_{i+1}$, the
  difference between the second and the third is
  \[
  \chi^{-1}[Q_{i-2}] - \chi^{-1}[Q'_i]=\chi^{-1}[Q_{i+1+2}] - \chi^{-1}[Q'_{i+1-1}],
  \]
  which is the image of $\ga_{i+1}$. Finally computing on $E_{i-1}$
  the difference between the third and the first is the image of
  \[
  \chi^{-1}[Q_{i+2}] - \chi^{-1}[Q'_i]=\chi^{-1}[Q_{i-1-2}] - \chi^{-1}[Q'_{i-1+1}]=
  \chi^{-1}[Q_{i-1+2}] - \chi^{-1}[Q'_{i-1-1}],
  \]
  which is the image of $\ga_{i-1}$, and where the last equality is by passing
  from two branch points of a $g^1_2$ on an elliptic curve to the two
  complimentary ones.
  Hence, the image of
  $\gb_i+\ga_{i+1}+\ga_{i-1}$ under $\prod \Pic E_i\to\Pic C$ is trivial; since
  $\gb_i$ is a 2-torsion point, this statement is equivalent to the claim
  in the Proposition. It remains to show
  that the span of the $\ga_i$ is an isotropic group w.r.t. the Weil
  pairing, which is clear since they are supported on different components
  in the product $\Pi E_i$.
\end{proof}
\begin{ntt}\label{N:a_prime}
  We identify the $\ga_i$ with their images in $\Pic^0 C[2]$, and denote
  by $\ga_i'$ the unique element whose Weil pairing with $\ga_j$ is non-trivial
  if and only if $j=j'$. I.e.\
  \[
  \Pic^0 C[2]=\spn(\ga_0,\ldots,\ga_4)\oplus\spn(\ga'_0,\ldots\ga'_4),
  \]
  where the only non-trivial Weil pairings between the ten basis elements
  $\ga_i, \ga'$ are the
  five pairings $\langle\ga_i,\ga'_i\rangle=1$.
\end{ntt}
\begin{prd}[``easy to describe'' effective theta characteristics]\label{P:easy}
  Let $E_i, E_j$ be such that $\chi(\Pic^2 E_i)\cap\chi(\Pic^2 E_j)=[Q_l]$,
  then the $\chi$ pullback of $[Q_l]$ is the intersection point of
  $\Pic^2 E_i\cap \Pic^2 E_j$ with multiplicity $2$. Denote this intersection
  point by $P$, and its pullbacks to $E_i, E_j$ by $p_{ij}, p_{ji}$ respectively.
  Denote by $q_{ij}^k$ (where $k=1,\ldots, 4$) the points
  on $E_i=\Pic^1 E_i$ so that $2q_{ij}^k=p_{ij}$; similarly denote by $q_{ji}^{k'}$
  the points on $E_j=\Pic^1 E_j$ so that $2q_{ji}^{k'}=p_{ji}$.
  Then:
  \[
  2(\pi_i^{-1}(q_{ij}^k) + \pi_j^{-1}(q_{ji}^{k'}))=
  \pi_i^{-1}(\ol{p_{ij}})+\pi_j^{-1}(\ol{p_{ji}})=2\theta_{ij}=K_C
  \]
  Hence $(\pi_i^{-1}(q_{ij}^k) + \pi_j^{-1}(q_{ji}^{k'}))$ is an effective
  theta characteristic.
\end{prd}
\begin{cor}\label{C:distinct}
  The theta characteristics constructed above are odd and distinct
  (hence there are $10\times 2^4$ of them).
\end{cor}
\begin{proof}
  It is clear that for each $\{i,j\}$ we get $2^4$ different effective
  theta characteristics, and that they are disjoint from the ten known effective
  even theta characteristics, which are the only effective even theta
  characteristics. Hence these $16$ theta characteristics are all odd.
  It remains to prove that the ten $16$-tuples are disjoint.

  Note that there is no point $p\in C$ with two vanishing coordinates:
  Indeed if there was such a $p$ then, using the fact that each of the $Q_i$s
  and $Q'_i$s are each supported on three coordinates, all the coordinates of
  $p$ would have to vanish.

  Note also that of the four points on $C$ comprising the theta characteristic
  in Proposition-Definition~\ref{P:easy},
  the two points in the $\pi_j$ pre-image of
  $q_{ij}^k$ have vanishing $i$ coordinate,
  their $j$ coordinate are interchangeable by $x_j\mapsto -x_j$, and their
  other coordinates are identical. Which means that --- by the claim above ---
  the theta $16$ tuple we get for a given $\{i,j\}$ is disjoint from the $16$
  tuple we get from a different choice $\{i', j'\}$.
\end{proof}
\begin{prd}
  By construction, each of the ten $16$-tuples of theta characteristics above
  is a translate of one of the spaces we now define:
  \[\begin{aligned}
  V_{ij}&:=\Pic E_i[2]\oplus\Pic E_j[2]=
  \spn(\ga_i, \ga_j,\gb_i,\gb_j)\\
  =&\spn(\ga_i, \ga_{i-1}+\ga_{i+1},\ga_j, \ga_{j-1}+\ga_{j+1})
  \end{aligned}\]
  In the case where $|i-j|=1$, in which case w.l.o.g. $i+1=j$,
  the last span is
  $span(\ga_i, \ga_{i-1}, \ga_{i+1}, \ga_{i+2}).$
  In the second case w.l.o.g. $i+2=j$ in which case the span is
  $span(\ga_i, \ga_{i-1}+\ga_{i+1}, \ga_{i+2}, \ga_{i+1}+\ga_{i+3})$.
\end{prd}
\begin{cor}
  Recalling Notation \ref{N:a_prime}, and taking indices mod 5, the following
  identities hold in $J C[2]$:
  \[
  V_{ij}^\perp=\begin{cases}
  \spn\{\ga_1,\ldots\ga_5\}\oplus\spn(\ga'_{i+3})&\text{for }j=i+1\\
  \spn\{\ga_1,\ldots\ga_5\}\oplus\spn(\ga'_{i-1}+\ga'_{i+1}+\ga'_{i+3})&\text{for }j=i+2,
  \end{cases}
  \]
  where orthogonality is w.r.t. the Weil pairing.
\end{cor}
\begin{proof}
  The orthogonality to $\spn\{\ga_1,\ldots\ga_5\}$ is immediate. The
  orthogonality of the basis elements to the last element in the respective
  bases is by direct verification.
\end{proof}
\begin{ntt}\label{N:eta}
  For each $i,j$ denote by $\eta_{ij}$ the unique non-trivial element in the
  projection of $V_{ij}^\perp$ on the second component in the direct sum
  decomposition
  \[
  J C[2]=\spn\{\ga_i\}_{i=1}^5\oplus\spn\{\ga'_i\}_{i=1}^5.
  \]
\end{ntt}
%
\section{Subsets of Steiner systems on Humbert curves}\label{S:steiner}
%
We refer the reader to Chapter 5 of \cite{D} for background material
on theta characteristics. Throughout this section we continue to work with a
Humbert curve $C$.
\begin{ntt}
  Denote by $O_{160}$ the set of 160 odd theta characteristics we found in
  Proposition-Definition~\ref{P:easy}.
  Denote by $O_{ij}$ the unique translate of $V_{ij}$ into $O_{160}$.
  We will call $\{\{\gth, \gth+\ga\}|\gth\in\gS_\ga\}$ the {\em set of pairs of the Steiner system $\gS_\ga$}.
\end{ntt}
The main objective of this section is to prove the following:
\begin{prp}\label{P:all2K}
  Assume that $A\subset{O_{160}\choose 4}$ is such that
  \[
  \forall a\in \binom{O_{160}}{4}\sm A:\sum_{\gth\in a}\gth\neq 2K_C.
  \]
  Let $R$ be the
  symmetric relation on ${O_{160}\choose 2}$ defined by $\gth_1R\gth_2$ if
  $\gth_1\cup\gth_2\in A$. Further assume that $R$ is transitive and that
  $R$ partitions $O_{160}\choose 2$ into $510$ equivalent classes, then these
  equivalent classes are the intersections of $O_{160}\choose 2$
  with the sets of pairs of the Steiner systems of $C$.
\end{prp}
For ``auxiliary'' objectives of this section, and their motivation, see the
discussion at the end of this section.

Recall (see e.g.\ \cite{D} 5.1 for a detailed overview) that for any curve $C$,
the set $J C[2]$ is a
$2g$-dimensional symplectic space over $\BF_2$ w.r.t. the Weil pairing, and that
the set of theta characteristics of $C$, together with
$\gth\mapsto h^0(\gth)\mod 2$, and with the addition action of $J C[2]$ on theta
characteristics, is a corresponding affine symplectic space. Further recall that
a subspace of a linear symplectic space $L\subset V$ is called {\em isotropic}
if the symplectic pairing is trivial on it, and that in this case the symplectic
pairing on $L$ induces a symplectic pairing on $L^\perp/L$.
\begin{prp}\label{P:unique_translate}
  Let $L$ be a $(g-1)$-dimensional isotropic subspace in a $2g$ dimensional
  linear symplectic
  space over $\BF_2$, then $L$ has a unique translation in the odd points of
  a $2g$ affine symplectic spaces.
\end{prp}
\begin{proof}
  Since $L$ is isotropic, it induces a symplectic pairing on
  the 2-dimensional (over $\BF_2$) symplectic space $L^\perp /L$.
  Moreover, there is a one to to one correspondence between such translated
  spaces and the odd
  points of an affine space over the 2-dimensional space $L^\perp /L$, of which
  there is exactly one.
\end{proof}
\begin{prp}\label{P:oij_int}
  Let $L$ be a maximal isotropic subspace of a $2g$ linear symplectic space
  over $\BF_2$; let $L_1, L_2$ be two distinct $g-1$ dimensional subspaces
  of $L$, let $l_1,l_2$ be the unique non-trivial elements in the images of
  $L_1,L_2$ in the 4 dimensional
  symplectic space $(L_1\cap L_2)^\perp / (L_1\cap L_2)$, and let $l'_1,l'_2$
  be so that $l_1'\cdot l_2=0,l_1'\cdot l_1=1,l_2'\cdot l_2=1,l_2'\cdot l_1=0$,
  where $\cdot$ is the symplectic pairing induced from the symplectic pairing on
  the original $2g$ dimensional space.
  Finally, let $O_1, O_2$ be the unique translates (per
  Proposition~\ref{P:unique_translate}) of $L_1, L_2$
  respectively into some affine $2g$ dimensional symplectic space over $\BF_2$.
  Then $\{x-y|x\in O_1, y\in O_2\}$ is the translate of $L$ by $l'_1+l'_2$.
  Moreover, each such difference is realized an equal number of times.
\end{prp}
\begin{proof}
  We endow the space $(L_1\cap L_2)^\perp / (L_1\cap L_2)$ with an affine
  symplectic structure given by $\gep_1l_1+\gep'_1l'_1+\gep_2l_2+\gep'_2l'_2\mapsto \gep_1\gep_1'+\gep_2\gep'_2$. This
  induces affine structures on $L_1^\perp/L_1$ and $L_2^\perp/L_2$.
  The pre-image of the unique odd point
  of the 2 dimensional symplectic space $L_1^\perp/L_1$ (resp. $L_2^\perp/L_2$)
  in the 4 dimensional symplectic space $(L_1\cap L_2)^\perp / (L_1\cap L_2)$
  under the natural projection is the pair $\{l_1+l_1'+l_2, l_1+l_1'\}$ (resp.
  $\{l_2+l_2'+l_1', l_2+l_2'\}$).
  Hence, the set of differences of these two
  pairs (here we compute the differences between points in
  an affine space, which is a linear space torsor, yielding a point
  in the linear space)
  is the translate by $l'_1+l'_2$ of the set $\{0, l_1, l_2, l_1+l_2\}$,
  which is the image of $L$ in the quotient space
  $(L_1\cap L_2)^\perp / (L_1\cap L_2)$ under the quotient map.
  Since all affine symplectic spaces of a given dimension are isomorphic, the
  claim about the identity of the set of differences
  follows. The claim about the number of representative for each
  difference being equal follows from the fact that the computations
  above were done in a quotient vector space.
\end{proof}
\begin{prp}\label{P:different_ij}
  For each non zero element
  $\mu\in\spn\{\ga'_i\}_{i=0}^4\cap(\sum_{i=0}^4\ga_i)^\perp$
  there are exactly three pairs of indices $\{i,j\}\neq\{i', j'\}$ so
  that the difference set $O_{ij}-O_{i'j'}$ equals to
  $\spn\{\ga_i\}_{i=0}^4+\mu$, where each element appears 8 times for this
  difference. Moreover, this description covers all the possible
  indices choices $\{i,j\}\neq\{i', j'\}$.
\end{prp}
\begin{proof}
  We work with indices mod 5.
  We start by naming the elements in $\spn\{\ga'_i\}_{i=1}^5\cap(\sum_{i=0}^4\ga_i)^\perp\sm\{0\}$:
  \[
  \gc_k^I:=\ga'_k+\ga'_{k+1},\quad
  \gc_k^{II}:=\ga'_k+\ga'_{k+2},\quad
  \gc_k^{III}:=(\sum_m\ga'_m)-\ga'_k.
  \]
  By Proposition~\ref{P:oij_int} the difference
  set $O_{ij}-O_{i'j'}$ equals to the translate of $\spn\{\ga_i\}_{i=1}^5$ by
  $\eta_{ij}+\eta_{i'j'}$, where each element in the translated span appears 8
  times as a difference.
  Recalling the definition of the $\eta_{ij}$s in Notation \ref{N:eta},
  the sum $\eta_{ij}+\eta_{i'j'}$ equals:
  \begin{itemize}
  \item if $j=i+1, j'=i'+1$ (and w.l.o.g. $i'-i\in\{1,2\}\mod 5$),
    \[
    \ga'_{i+3}+\ga'_{i'+3}=\begin{cases}\gc^I_{i+3}&\text{for }i'=i+1\\
    \gc^{II}_{i+3}&\text{for }i'=i+2
    \end{cases},
    \]
  \item if $j=i+2, j'=i'+2$, (and w.l.o.g. $i'-i\in\{1,2\}\mod 5$),
    \[
      \ga'_{i-1}+\ga'_{i+1}+\ga'_{i+3}+\ga'_{i'-1}+\ga'_{i'+1}+\ga'_{i'+3}
      =\begin{cases}
      \gc^{III}_{i-1}&\text{for }i'=i+1,\\
      \gc^I_{i-1}&\text{for }i'=i+2,
      \end{cases}
    \]
  \item if $j=i+1, j'=i'+2$,
    \[
    \ga'_{i+3}+\ga'_{i'-1}+\ga'_{i'+1}+\ga'_{i'+3}=\begin{cases}
    \gc^{II}_{i'-1}&\text{for }i'=i\\
    \gc^{III}_{i'}&\text{for }i'=i+1\\
    \gc^I_{i'+3}&\text{for }i'=i+2\\
    \gc^{III}_{i'+2}&\text{for }i'=i+3\\
    \gc^{II}_{i'+1}&\text{for }i'=i+4.
    \end{cases}
    \]
  \end{itemize}
  The claim about
  each $\mu$ appearing $3$ times follows by counting. The claim about the
  multiplicity was already addressed, and the last claim
  is clear since $\binom{10}{2}=45=3\cdot(16-1)$.
\end{proof}
\begin{prp}\label{P:same_ij}
  The difference set $O_{ij}-O_{ij}$ is $V_{ij}$, where each
  difference is represented $\#V_{ij}$ times; hence the difference set in
  $O_{ij}\choose 2$ is $V_{ij}\sm\{0\}$, where each difference is represented
  $\# V_{ij}/2=16/2=8$ times.
\end{prp}
\begin{proof}
  The first part holds since  $O_{ij}-O_{ij}=V_{ij}-V_{ij}$, and $V_{ij}$ is a
  vector space; the second part follows from symmetry.
\end{proof}
\begin{prp}\label{P:same_ij2}
  For any $\mu\in\spn\{\ga_i\}_{i=1}^5\sm\{0,\sum_{i=0}^4\ga_i\}$,
  $\mu$ appears as a difference in one of the ten sets
  $(O_{ij}-O_{ij})\sm\{0\}$ precisely 48 times if $\mu$ is of the form
  $\ga_k,\ga_k+\ga_{k+2}$ or $\ga_{k+2}+\ga_{k+3}+\ga_{k+4}$, and precisely 32
  times otherwise, where we take indices modulo 5.
\end{prp}
\begin{proof}
  We start by naming the non $\ga_i$ elements in
  $\spn\{\ga_i\}_{i=1}^5\sm\{0,\sum_{i=0}^4\ga_i\}$:
  \[
  \begin{aligned}
  \ol{\gc}_k^I:=&\ga_k+\ga_{k+1},\quad
  \ol{\gc}_k^{II}:=\ga_k+\ga_{k+2},\quad
  \ol{\gc}_k^{III}:=(\sum_i\ga_i)-\ga_k,\\
  \ol{\gc}_k^{IV}:=&(\sum_i\ga_i)-\ol{\gc}_k^I,\quad
  \ol{\gc}_k^V:=(\sum_i\ga_i)-\ol{\gc}_k^{II}.
  \end{aligned}
  \]
  Then the following equalities hold by direct verification:
  {\small
  \[
  \begin{aligned}
    V_{i,i+1}=&\spn\{\ga_i\}_{i=1}^5\cap{\ga'_{i+3}}^\perp\\
    =&\{0,\ol{\gc}_i^I,\ol{\gc}_{i+1}^I,\ol{\gc}_{i+4}^I,
  \ol{\gc}_i^{II},  \ol{\gc}_{i+2}^{II},  \ol{\gc}_{i+4}^{II},
  \ol{\gc}_{i+3}^{III},\ol{\gc}_{i+3}^{IV},\ol{\gc}_{i+2}^{IV},
  \ol{\gc}_{i+3}^V,\ol{\gc}_{i+1}^V\}\cup
  \{\ga_j\}_{j\neq i+3},\\
  V_{i, i+2}=&\spn\{\ga_i\}_{i=1}^5\cap(\ga'_{i-1}+\ga'_{i+1}+\ga'_{i+3})^\perp\\
  =&\{0,\ga_i,\ga_{i+2},\ol{\gc}_{i+3}^I,
  \ol{\gc}_{i-1}^{II},  \ol{\gc}_i^{II},  \ol{\gc}_{i+1}^{II},
  \ol{\gc}_{i+1}^{III}, \ol{\gc}_{i-1}^{III}, \ol{\gc}_{i+3}^{III},
  \ol{\gc}_{i+4}^V,\ol{\gc}_{i+2}^V\}\cup\{\ol{\gc}_j^{IV}\}_{j\neq i+3}.
  \end{aligned}
  \]
  }
  Observe that the $\ga_\bullet,\ol{\gc}^{II}_\bullet,\ol{\gc}^{IV}_\bullet$s
  appears six times above, whereas the $\ol{\gc}^I_\bullet,\ol{\gc}^{III}_\bullet,\ol{\gc}^{V}_\bullet$s appear four times. This means that when we consider the
  orbit of the above sets under coordinate rotation, we will have
  $30$ occurrences of elements of the first three types, and $20$ of the second.
  However, as each class of element is also of size five, and as rotations
  act transitively on each class, we see that each of the
  elements in the first classes is represented six times, and each of the
  elements in the later classes is represented four times. The result now
  follows from Proposition~\ref{P:same_ij}.
\end{proof}
\begin{rmr}\label{R:gamma_orbits}
  The orbit of each of the $\mu$'s from Proposition~\ref{P:same_ij2} under
  either the automorphism group
  of $C$, or its index $2$ subgroup generated by the rotations and sign
  changes on the coordinates, is of size 5.
\end{rmr}
\begin{proof}
  These groups act on the $\ga_k$ via the action of dihedral group and the
  cyclic group on the indices, respectively.
\end{proof}
\begin{proof}[Proof of Proposition~\ref{P:all2K}]
  We start by noting that
  the set of differences between pairs of distinct elements of $O_{160}$
  is $(\sum_i\ga_i)^\perp\sm\{0,\sum_i\ga_i\}$:
  Indeed, in Propositions~\ref{P:different_ij}, \ref{P:same_ij} we showed that
  \[
  \cup_{\{i,j\}\neq\{i',j'\}}O_{i,j}-O_{i',j'}=\left(\sum_i\ga_i\right)^\perp\sm\spn(\{\ga_i\}_i),
  \]
  whereas in \ref{P:same_ij2} we showed that
  \[
  \cup_{\{i,j\}}O_{i,j}-O_{i,j}=\spn(\{\ga_i\}_i)\sm\{(\sum_i\ga_i),0\}.
  \]
  By the definition of $A$, no members of two equivalence classes of $R$
  are in the same set of pairs of some Steiner system. Moreover, since
  $O_{160}\choose 2$
  has representatives in sets of pairs of $510$ Steiner systems, none of the
  equivalent classes can break into (sets of pairs of) Steiner systems.
\end{proof}
  Sadly, this is as far as I managed to proceed with the
  non-numerical computation of the Steiner systems. It is important to
  note that the requirement for the certification of the numerical part
  is merely the number of partial Steiner systems represented in
  $\binom{O_{160}}{2}$, which is somewhat weaker
  than the enumerative parts of Propositions~\ref{P:different_ij} and
  \ref{P:same_ij}; we use the later ones, as well as
  Remark~\ref{R:gamma_orbits}, for verification purposes (see \ref{D:tests}).
%
\section{The Wiman curve}\label{S:rep}
%
\begin{prd}[Wiman curve $\bW=W^{160}$]\label{PD:wiman}
In \cite{W}, Wiman discovered a genus $5$ curve
(denoted in \cite{E} by $W^{160}$) whose canonical model is the intersection
of the nulls of the quadric forms:
\[
  Q_A:=\sum_{i=0}^4 x_i^2, \quad Q_B:=\sum_{i=0}^4\zeta^ix_i^2,\quad Q_C:=\sum_{i=0}^4\zeta^{-i}x_i^2,
\]
where $\zeta^5=1,\zeta\neq 1$. The automorphism group of this curve is
the semi-direct product of $(\BZ/2)^5/(\BZ/2)$, acting by multiplying by $-1$ each
of the coordinates, and the dihedral group $D_5$ acting on the
coordinates, generated by
$(x_0,x_1,x_2,x_3,x_4)$, and $(x_0,x_3)(x_1,x_2)$. From now on we assume that
the canonical system of the Wiman curve is endowed with these coordinates.
\end{prd}
As explained in the introduction, the object of this section is two-fold:
we explicitly compute the coordinates of $O_{ij}$ for the curve $W$, and
we decompose the quadrics forms on the dual canonical system of $W$ to
irreducible representations of the automorphism group.
\begin{prd}[The even effective theta characteristics of $\bW$]
  Denote by $\phi$ the ``golden ratio'' $(1+\sqrt{5})/2$. W.l.o.g.\:\:we can
  assume that $\zeta=\exp(2\pi i / 5)$, in which case
  \[
  \zeta+\ol{\zeta}=\phi - 1,\quad \zeta^2+\ol{\zeta}^2=-\phi.
  \]
  We now compute
  \[\begin{aligned}
    Q_B +& Q_C = 2x_0^2 + (\zeta+\ol{\zeta})(x_1^2 + x_4^2) + (\zeta^2+\ol{\zeta}^2)(x_2^2 + x_3^2)\Rightarrow\\
    Q_0:=&\frac{\phi}{\phi^2+1}\left(Q_B+Q_C+\phi Q_A\right)=
    \frac{\phi}{\phi^2+1}\left((2+\phi)x_0^2+(2\phi-1)(x_1^2+x_4^2)\right)\\
    =&
    \frac{\phi}{\phi^2+1}\left((\phi^2+1)x_0^2+(\phi + 1/\phi)(x_1^2+x_4^2)\right)=
    \phi x_0^2+(x_1^2+x_4^2),\\
    Q_0':=&\frac{1}{\phi-3}(Q_B+Q_C-(\phi-1)Q_A)\\
    =&\frac{1}{\phi-3}\left((2-(\phi-1))x_0^2+(-\phi-(\phi-1))(x_2^2+x_3^2)\right)
    =-x_0^2+\phi(x_2^2+x_3^2),
  \end{aligned}\]
  where the last equality follows since
  $\phi-3=\frac{1}{\phi}-2=\frac{1-2\phi}{\phi}$.
  One gets the other $8$ even effective theta characteristics $Q_i, Q_i'$
  by applying the cyclic group (which is a subgroup of the symmetry group
  of the curve) on the coordinates.
\end{prd}
\begin{prp}\label{P:wiman_odd}
  Identify the curve $\bW$ with its canonical image, and denote
  $a:=1/\sqrt{\phi},\ \ii:=\sqrt{-1}$. Then the effective theta
  characteristics from \ref{P:easy}
  are the orbit under the cyclic action $\BZ/5\BZ$ on the coordinates, (as a
  subgroup of the automorphism group of the curve) of the $32$ divisors:
  {\small
  \[\begin{aligned}
  (\ii;\gep_1a;0;\ii a,\gep_2)+(\ii;\gep_1a;0;-\ii a,\gep_2)+
  (\ii;\gep_3;\ii a;0;\gep_4a)+(\ii;\gep_3;-\ii a,0;\gep_4a),\\
  (0;ia;\gep_1;\ii;\gep_2a)+(0;\ii a;\gep_1;-\ii;\gep_2a)+
  (\ii;\gep_3;\ii a;0;\gep_4a)+(-\ii;\gep_3;\ii a;0;\gep_4a),
  \end{aligned}\]
  }
where $\gep_1,\gep_2,\gep_3,\gep_4\in\{-1,1\}$ are sign choices.
\end{prp}
\begin{proof}
  We find the first type of point by working (in the terminology of
  Proposition~\ref{P:easy}) with $i,j=2,3$; i.e.\ the construction
  corresponding
  to starting with the even theta characteristic corresponding to $Q_0$.
Recall (see Proposition-Definition~\ref{p:elliptic_cone}) that $E_3$ is the base of the cone $Z(Q_0)\cap Z(Q_4')$,
whereas $E_2$ is the base of the cone $Z(Q_0)\cap Z(Q_1')$. We proceed
to find $q_{32}^k$ (and later $q_{23}^k$).
We start by finding the ramification points of the map from
$E_3$ to the base of the cone $Z(Q_0)$:

Assuming $x_0=0$ we get
\[
Q_0=0\text{ and }x_0=0\quad\Leftrightarrow\quad x_4^2 =-x_1^2;
\]
substituting in the equation of $Q_4'$ we get $0=\phi x_1^2 + \phi x_2^2 + x_1^2$,
and the ramification points of the map from $E_3$ to the base of the cone
$Z(Q_0)$ are the locus of the
trivial $x_2$ discriminant; however the discriminant for $x_2$ is $0$ if and
only if $x_1=x_2=0$, which is impossible as at least one of the coordinates has
to be non zero.

Hence $x_0\neq 0$. Since we work with homogeneous equations, we may choose
$x_0\neq 0$ arbitrarily. Setting $x_0^2=-1/\phi$ and substituting in $Q_0$ we
get:
\[
x_4^2=1-x_1^2,\quad\Rightarrow\quad
(\text{substituting $x_4$ in }Q_4'):\quad0=\phi x_1^2 + \phi x_2^2 - (1-x_1^2),
\]
where the discriminant for $x_2$ is $0$ when $x_2=0$, and then
\[
1=\phi^2 x_1^2\quad\Rightarrow\quad x_1=\pm1/\phi\quad\Rightarrow\quad
x_4^2=1-1/\phi^2=1/\phi.
\]
To find the pullbacks of these points to $C$ we solve $x_3$ in $Q_1'$, giving us
\[
 0=\phi x_3^2+\phi(1/\phi)-1/\phi^2
 =\phi x_3^2+1-1/\phi^2
 =\phi x_3^2+1/\phi\Leftrightarrow x_3=\pm \ii/\phi.
 \]
 The coordinates written in the statement of the Proposition are projectively
 equivalent.
 To find $q_{23}^k$ we apply the symmetry switching between $Q_4'$ and $Q_0'$;
 i.e.\ switching  $x_1,x_4$ and $x_2,x_3$.

 As for the other type we consider $i,j=3,0$ (i.e.\ the construction
 corresponding to starting with the theta characteristic corresponding to
 $Q_4'$):
We now find the ramification points of the map from $E_3$ to the base of the
cone $Z(Q_4')$.
As before, to solve $Q_4'=0$ there are two cases:
similarly to the case above, the case where $x_4=0$ does not give
a ramification point: indeed assuming $x_4=0$ we get $0=x_1^2+0+\phi x_0^2$,
and the $x_0$ discriminant is $0$ if and only if both are $0$.
In the other case we set $x_4^2=\phi$, and then $x_1^2=1-x_2^2$, and we
have to find when the $x_0$ discriminant of
$0=1-x_2^2+\phi+x_0^2$ is $0$, which happens when
\[
x_0=0, x_2^2=\phi^2\Rightarrow x_1^2=-\phi.
\]
To find the pullbacks we solve $Q_3$:
$0=\phi^2+\phi+\phi x_3^2\Rightarrow x_3^2=-\phi^2$.
To find $q_{03}^{k'}$ we apply the symmetry switching $Q_0$ and $Q_3$; i.e.\ we
switch between $x_0,x_3$, and $x_1,x_2$.
\end{proof}
\begin{prp}\label{P:irrep_decomp}
  Taking indices mod 5,
  The decomposition of $\Sym^2H^0(K_\bW)$ into irreducible representations of $(\BZ/2)^5\rtimes D_5$ (where the $k$th copy of $\BZ/2$ acts by $x_k\mapsto-x_k$) is
  \[
  \begin{aligned}
    \ &\spn(\{x_jx_{j+1}\}_{j=1}^5), \quad \spn(\{x_jx_{j + 2}\}_{j=1}^5),\quad\spn(\sum_jx_j^2),\\
    \ &\spn(\sum_je^{\ii\pi j/5}x_j^2,\sum_je^{-\ii\pi j/5}x_j^2),\quad
    \spn(\sum_je^{2\ii\pi j/5}x_j^2,\sum_je^{-2\ii\pi j5}x_j^2).
  \end{aligned}
  \]
  Moreover, all these representations are not isomorphic; not even as
  projective representations.
\end{prp}
\begin{proof}
  Taking indices modulo 5, $k\in\{1,2\}$ and given a non-trivial vector $v$ in
  $\spn\{x_jx_{j+k}\}_{j=1}^5$, assume that the coefficient in of $x_mx_{m+k}$ in $v$
  is non zero, then the sum of $v$ with its image under
  $x_{m-1}\mapsto-x_{m-1},x_{m+k+2}\mapsto-x_{m+k+2}$ is twice this coefficient times
  $x_mx_{m+k}$, and the orbit of this element under the group spans
  $\spn\{x_jx_{j+k}\}_{j=1}^5$.

  To show that these two irreducible representations (for $k=1,2$) are not
  isomorphic, we consider the decomposition to $-1$ and $1$ eigen-spaces under
  the
  action of $x_0\mapsto-x_0, x_2\mapsto-x_2$. Indeed we see that these
  decomposition are
  \[
  \begin{aligned}
  \spn\{x_jx_{j+1}\}_{j=1}^5=&\ \spn\{x_{j-1}x_j\}_{j=0}^3\bigoplus\spn x_3x_4,\\
  \spn\{x_jx_{j+2}\}_{j=1}^5=&\ \spn\{x_3x_0,x_2x_4\}\bigoplus\spn\{x_4x_1,x_0x_2,x_1x_3\},
  \end{aligned}
  \]
  i.e.\ they are not isomorphic.

  Finally, the group $(\BZ/2)^5$ acts trivially on $\spn(\{x_j^2\}_{j=0}^4)$,
  so for this subspace we only have to consider the classical representation
  theory of dihedral groups.
\end{proof}
\begin{cor}\label{C:irrep_decomp}
  The projectivization of the span of the 3rd and 4th representations in
  Proposition~\ref{P:irrep_decomp} is $I_2(\bW)$.
\end{cor}
\begin{proof}
  Direct verification.
\end{proof}
%
\section{Certifying a divisor is not two-canonical}\label{S:approx}
%
The object of this section is to provide a certification algorithm
for a divisor in $\Div^{16}(W)$ to be non two-canonical, with some assumption
on the maximal multiplicity it incurs (specifically, we deal with all cases
with multiplicity at most 3).
By ``describing the algorithm'' we
mean a combination of three things: we review all the
approximate computations need to be carried out, give bounds on how accurate
they have to be, and finally describe how to perform the computations.

From here on we identify all the fibers of the tangent bundle
of the linear space $H^0(K_\bW)^*$ with the tangent space at $0$, which is
itself isomorphic to $H^0(K_\bW)^*$.

Throughout this section we don't work with points in the dual canonical
system $|K_\bW|^*$, but rather with their pullbacks under the projectivization
map $H^0(K_\bW)^*\sm\{0\}\to|K_\bW|^*$, or with representatives of these points;
i.e.\ points along the fiber of the projectivization map.

A final point to keep in mind is that the object of this section is to build
an algorithm to solve a multi-linear algebra problem, and that
the prevalent linear algebra libraries which work over the complex numbers
work with Hermitian dot products; hence the relevant statements below
are sometimes made in this language. For the same reason we sometimes use
complex conjugation, denoted by $x\mapsto\ol{x}$.
\begin{lma}\label{L:tangent}
  Let $p\in H^0(K_\bW)^*$ be a representative of some point $\hat{p}\in W$, with a
  vanishing coordinate $j$. Then
  the pullback of the tangent $T_{\hat{p}}W$ to $H^0(K_\bW)^*$ is linearly spanned by
  two vectors:
  \begin{itemize}
  \item The vector $\ol{p}$, which spans the fiber of the
    projectivization over $\hat{p}$,
  \item and the vector whose only non zero coordinate is $j$, which is
    orthogonal to it (in either Hermitian or non-Hermitian dot product).
  \end{itemize}
\end{lma}
\begin{proof}
  The claim about the first vector is immediate.
  As for the second, note that
  since all the $Q_k$ are of the form $\sum a_ix_i^2$ --- i.e.\ diagonalized ---
  their gradients are of the form
  \[
  2(a_0x_0, a_1x_1, a_2x_2, a_3x_3, a_4x_4);
  \] specifically, their $j$ coordinate is trivial when evaluated at $p$.
  The orthogonality claim is also immediate.
\end{proof}
\begin{prd}\label{PD:Hessian}
  Let $p,\hat{p},j$ as in Lemma~\ref{L:tangent} above.
  Let $\{u_n\}_{n=0}^2$ be a basis for $I_2(W)$, and
  let $f$ be a homogeneous polynomial on $H^0(K_\bW)^*$ which attains a double
  zero at $\hat{p}\in W$; i.e.\
  \[
  f(p) = 0,\quad\nabla f|_p=\sum_n \lambda_n\nabla u_n|_p
  \text{ for some }\lambda_n\text{s}.
  \]
  Define the Lagrangian $\cL(x_k, \lambda_n):=f-\sum_n\lambda_nu_n$,
  then $f$ has a triple zero on $W$ at $\hat{p}$ if and only if the $j,j$ entry
  of the partial Hessian (i.e.\ the Hessian w.r.t. the coordinates $x_i$)
  $\Hess_x\cL$ is $0$.
\end{prd}
\begin{proof}
  We use a very slight modification of the bordered Hessian argument
  used in constrained optimization:
  Let $r(t)$ be an analytic function from the unit disc in $\BC$ to
  $H^0(K_\bW)^*$, so that
  $r(0)=p$, and $r(t)$ projects to $W$ under the projectivization map.
  Denote the gradient of $r(t)$ by $v$, then applying the chain rule we have
  \[
  \frac{d}{dt^2}f(r(t))=
  v^T\cdot \Hess\,f\cdot v +\nabla f\cdot r''(0),
  \]
  where $-^T$ denotes the transpose, and $\cdot$ denotes the coordinate-wise
  product (and not the Hermitian bi-linear form, i.e.\ there is no conjugation
  involved). Likewise we get:
  \[
  \frac{d}{dt^2}u_n(r(t))=
  v^T\cdot \Hess\,u_n\cdot v +\nabla u_n\cdot r''(0).
  \]
  Since $\nabla f =\sum_n\lambda_n \nabla u_n$ at $p$, we now have
  \[
  \frac{d}{dt^2}f(r(t))|_{t=0}=
  v^T\cdot\left(\Hess\,f-\sum_n\lambda_n\Hess\,u_n\right)\Big|_{t=0}\cdot v=
  v^T\cdot\Hess_x\cL|_{t=0}\cdot v.
  \]
  Hence the function $f$ has a triple $0$ at $\hat{p}$ if and only if the right
  hand side of the above equation is $0$, where $v$ is a representative in
  the fiber of the
  pullback of $T_{\hat{p}}W$ to $H^0(K_W)$,
  which is orthogonal to the fiber of the projectivization
  map.
  By Lemma~\ref{L:tangent}, we are done.
\end{proof}
\begin{dsc}[Efficient computation of the $j,j$ coefficient of the Hessian, given $\nabla f$]\label{D:half_grad_coef}
  In the notations of Proposition~\ref{PD:Hessian},
  taking indices modulo 5, we need to compute the value of the $j,j$ entry of
  $\sum_n\lambda_n\Hess\,u_n|_p$:

  The $u_n$s we choose are simply $Q_{j-1}, Q_j, Q_{j+1}$. We now have to
  solve the equation
  \[ \nabla f=\lambda_0 \nabla Q_{j-1}+\lambda_1 \nabla Q_j + \lambda_2
  \nabla Q_{j+1}.
  \]
  Recalling that $Q_k=x_{k-1}^2+\phi x_k^2+x_{k+1}^2$ (for $k=0,\ldots, 5$),
  we see that only $Q_{j-1}$ (resp. $Q_{j+1}$) is contributing to the $j-2$
  (resp $j+2$) entry of the gradient. Hence, $\lambda_0$ (resp.  $\lambda_2$) is
  the $j-2$ (resp. $j + 2$) entry of $\nabla f$ divided by
  $2$ times the $j-2$ (resp $j + 2$) entry of $p$.

  In a similar manner, we see that $\lambda_1$ can be expressed either as the
  $j-1$ coefficient
  of $\nabla f$ divided by $2$ times the $j-1$ coefficient of $p$, minus
  $\phi\lambda_0$, or as the $j+1$ coefficient
  of $\nabla f$ divided by the $2$ times the $j+1$ coefficient of $p$, minus
  $\phi\lambda_2$.
\end{dsc}
\begin{prp}\label{P:whats2K}
  Endow $\Sym^2H^0(K_\bW)$ with the complex vector norm coming from taking the
  monomial basis on the $x_i$'s. Let $q\in\Sym^2H^0(K_\bW)$ be a quadric which
  does not vanish identically on $W$. Let $D\in\Div^{16}(W)$ be the intersection
  divisor of $Z(q)$ and the canonical model of $W$. Assume $p\in \BA_5\sm\{0\}$
  is a representative of some point in $D\subset |K_\bW|^*\cong\BP^4$, and let
  $p',q'$ be
  such that $||p-p'||<\gep, \quad ||q-q'||<\gd$, then the following properties hold:
  \begin{enumerate}
  \item\label{single} $||q'(p')||\leq15\left((||p||+\gep)^2(||q||+\gd)-||p||^2||q||\right)$.
  \item\label{double} If $p$ is a double point of $D$, then the norm of the
    unitary projection of $\nabla q'(p')$ on the unitary complement of
    $\spn\{\nabla Q_j(p)\}_j$ is bounded by
    $25\left((||p||+\gep)(||q||+\gd)-||p||\cdot||q||\right)$.
  \item\label{triple} If $p$ is a triple point of $D$, and
    $j,\lambda_n$ are as in \ref{D:half_grad_coef}, and
    let $\lambda_n'$ be $\gc$ approximations of the $Q_{j-1},Q_j, Q_{j+1}$
    coefficients of the projection of $q'(p')$ on $\spn\{\nabla Q_k(p)\}$,
    then the norm of element at the diagonal $j, j$ entry of
    $\Hess(q')-\sum_{n=0}^2\lambda_n\Hess(Q_{j-1+n})$ is bounded by
    \[
    2\gd+6\phi\gc+50\left((||p||+\gep)(||q||+\gd)-||p||\cdot||q||\right).
    \]
  \end{enumerate}
\end{prp}
\begin{proof}
  Denote the coordinates of $q,q'$ by $q_{i k},q'_{i k}$, and the coordinates
  of $p,p'$ by $p_i,p'_i$.
  As for the claim \ref{single}, observe that
  \[q'(p')-q(p)=\sum_{i\leq k}(q'_{i k}p'_ip'_k-q_{i k}p_ip_k),\]
  and apply the triangle inequality.

  Denote the unitary projection on the unitary complement of
  $\spn\{\nabla Q_j(p)\}_j$ by $P^\perp$. By Lagrange multipliers,
  $P^\perp\nabla q = 0$. Moreover, the projection of the $i$th component
  of $\nabla q'(p')-\nabla q(p)$ is given by
  \[
  P^\perp\left(\sum_k q'_{i k}p'_k-\sum_k q_{i k}p_k\right)=P^\perp\left(\sum_k(q'_{i k}p'_k-q_{i k}p_k)\right).
  \]
  As before, we now apply the triangle inequality; summing the results over
  $i=0,\ldots,4$, this proves claim \ref{double}.

  Finally, denote by $\{\lambda_n''\}_{n=0}^2$ the $Q_{j-1},Q_j, Q_{j+1}$
  coefficients of $P^\perp\sum_kq'_{i k}p'_k$.

  we have to bound the norm of the $j,j$ entry of the $5\times 5$ matrix
  \[
  \begin{aligned}
  \ &(\Hess(q')-\sum_{n=0}^2\lambda'_n\Hess(Q_{j-1+n}))\\
  =&(\Hess(q')-\sum_{n=0}^2\lambda'_n\Hess(Q_{j-1+n}))-
  (\Hess(q)-\sum_{n=0}^2\lambda_n\Hess(Q_{j-1+n}))\\
  =& \Hess(q'-q)-\sum_{n=0}^2((\lambda'_n-\lambda''_n)+(\lambda''_n-\lambda_n))
  \Hess(Q_{j-1+n}).\\
  =& \Hess(q'-q)-\sum_{n=0}^2(\lambda'_n-\lambda''_n) \Hess(Q_{j-1+n})
  -\sum_{n=0}^2(\lambda''_n-\lambda_n)\Hess(Q_{j-1+n}).
  \end{aligned}
  \]
  We will now bound the norm of each of these three terms separately.
  As for the first term any entry of $\Hess(q'-q)$ is of
  norm $<2\gd$.
  Similarly, for the second term, we note that all entries of $Q_m$ are of norm
  $<\phi$, so any entry of $\Hess(Q_m)$ is of norm $<2\phi$.
  Finally, to bound the third term we note that
  \[
  ||\sum_{n=0}^2(\lambda''_n-\lambda_n)\nabla Q_{j-1+n}||<||\nabla q'(p')-\nabla q(p)||,
  \]
  and that $Q_m$ are diagonal quadric forms, so up to a factor of $2$ this is also
  the bound on the Hessian of the sum. Hence, any entry of
  the difference above is bounded by
  $2(\gd+3(\phi\gc+||\nabla q'(p')-\nabla q(p)||))$, and
  we have already bound the last term in the proof of claim~\ref{double}.
  Hence, by Proposition-Definition~\ref{PD:Hessian}, we have proved claim
  \ref{triple}
\end{proof}
Having answered the ``what we have to compute'' and ``how accurate do we have the compute'' questions, we now move to the ``how do we compute'' questions.
Specifically we have to show how we intend to compute the following:
\begin{itemize}
\item Find the quadrics vanishing on a given set of points in
  $W\subset\BP^4$,
\item of the quadrics above, which are non-trivial under the quotient by
  $I_2(W)$,
\item of the above, which are non-trivial under the projection to the
  direct sum of tangents space of $W$ at given points,
\item of the above, which quadrics satisfy the condition for triple points.
\end{itemize}
We address the first three issues in Proposition~\ref{P:single_svd},
Corollary~\ref{C:single_svd}, and Corollary~\ref{C:svd_double}
respectively. With the heavy lifting already done, in all of these cases
it is more complicated to build the computational gadgets then the prove they
are the correct ones -- this is usually immediate. As for the last issue,
by Proposition-Definition~\ref{PD:Hessian}, and since we are dealing
with Hessians of quadratic forms, it does not involve solving or setting up
any linear
system; merely plugging a few numbers we already computed in a known formula.
\begin{ntt}
  Since we denote coordinates of points by subscripts, we need an alternative
  notation for indices when using several points. Out alternative is
  $-_{[i]}$. Note that there is no ambiguity in the case of quadric forms
  as in that case coordinates are double indices, so as indexed quadrics are
  used throughout in this paper we will not clutter the notation.
\end{ntt}
\begin{prp}\label{P:single_svd}
  Assuming $\{p_{[i]}\}_{i=1}^n$ are representatives of points in $W\subset\BP^4$,
  let $M_{15}$ be the $n\times 15$ matrix whose rows are
  the degree $2$ monomials in the coordinates of the points.
  Let $M_{15}=U_{15}D_{15}V_{15}^\dagger$ be
  the singular value decomposition $M_{15}$, then the rows of $V_{15}^\dagger$
  corresponding to the singular value $0$ --- considered as quadrics over $\BP^4$
  in lexicographic basis --- span the quadrics which vanish on all the
  $p_{[i]}$s.
\end{prp}
\begin{proof}
the right kernel of $M_{15}$
is the set of enveloping quadrics of the points, expressed in
lexicographic basis.
\end{proof}
\begin{cor}\label{C:single_svd}
  Let $M$ be the projection of $M_{15}$ on the space orthogonal --- per
  Proposition~\ref{P:irrep_decomp} and Corollary~\ref{C:irrep_decomp} ---
  to $I_2(W)$, and let $M=UDV^\dagger$ be
  the singular value decomposition $M$, then the rows of $V^\dagger$
  corresponding to the singular value $0$ --- considered as quadrics over $\BP^4$
  in lexicographic basis --- span the two canonical quadrics which vanish
  on all the $p_{[i]}$s, and are orthogonal (in representation theoretic sense)
  to $I_2(W)$.
\end{cor}
\begin{proof}
  This follows from Corollary~\ref{C:irrep_decomp}.
\end{proof}
\begin{rmr}[effective computation of gradient of a quadric]\label{R:grad_comp}
  Recall that if $q=\sum_{k,k'}q_{kk'}x_kx_{k'}$ is a quadric form as above, then
  the $k$th coordinate of $\nabla q|_p$ is $2q_{kk}p_k+\sum_{k'\neq k}q_{kk'}p_k'$.
\end{rmr}
\begin{cor}\label{C:svd_double}
  Let $q=\sum_{i=1}^{k_1}\lambda_iq_i$ be a sum of quadrics represented by the
  kernel
  columns in Proposition~\ref{P:single_svd}, such that the double points of
  $Z(q)\cap W$ are supported on points in $\BP^4$ represented by
  $\{p_{[j]}\}_{j=1}^{n_2}\subset\BC^5$. Let $M_2$ be the $2n_2\times k_1$ block matrix,
  whose $2\times k_1$ blocks are the dot products of the $\nabla q_i|_{p_{[j]}}$
  --- computed as in Remark~\ref{R:grad_comp} above --- with
  the two vectors from Lemma~\ref{L:tangent}, evaluated at $p_{[j]}$, and let
  $M_2=U_2D_2V_2^\dagger$ be the singular value decomposition of $M_2$, then
  $M_2$ has a $0$ singular value, and $\{\lambda_i\}$ is in the span
  of the rows of $V^\dagger$ corresponding to the $0$ values.
\end{cor}
\begin{proof}
  This follows from part \ref{double} in Proposition~\ref{P:whats2K}, and from
  the computation of the tangent space in Lemma~\ref{L:tangent}.
\end{proof}
\begin{rmr}
  The converse of Proposition~\ref{P:whats2K} is not necessarily correct even
  when $D$ is the sum of four theta characteristics in $O_{160}\choose 4$ and for
  $\gep,\gd=0$. Namely, one cannot use the opposite inequalities to deduce that
  the sum of four theta characteristics in $O_{160}\choose 4$ {\em is} $2K_\bW$.
  The argument is simply insufficient if the four theta characteristics share
  a common point, and the ``tests'' for single, double, and triple multiplicity
  all pass.

  The usage of the results of Corollary~\ref{C:svd_double} as an ``input''
  for \ref{D:half_grad_coef} makes things even worse in this
  respect: it assumes a one dimensional space of solutions, which is not
  guaranteed.

  Luckily, the combination of Propositions~\ref{P:whats2K} and
  \ref{P:all2K} show that one can substitute exact computation in a specific
  case by approximate computation in many cases.
\end{rmr}

%
\section{Certified numeric proof of the main theorem}\label{S:prog}
%
\begin{dsc}[Floating point representation and operations accurate]\label{D:ieee754}
  We represent floating point numbers using IEEE-754 double accuracy which has
  a $52$
  bits mantissa; i.e.\ representation accuracy of $2^{-53}\sim1.11\cdot 10^{-16}$.
  IEEE-754 guarantees at most half a bit
  error (i.e.\ $2^{-53}$ multiplicative error) for each multiplication, addition,
  or subtraction.
\end{dsc}
\begin{dsc}[Bounding the accuracy of SVD]\label{D:svd_accuracy}
  All the SVD decomposition we use are accurate in the following sense:
  Let $A$ be a matrix {\em which occurs in our program} (described below), for
  which we compute the SVD (these matrices sizes are at most $12\times 16$ in
  Corollary~\ref{C:single_svd}, at most $8\times 16$ in
  Corollary~\ref{C:svd_double}, and at most
  $15\times 48$ in verifying Corollary~\ref{C:irrep_decomp}). Denote the output
  of the floating point SVD by $U,D,V$, then all the entries of
  $UU^\dagger-Id,VV^\dagger-Id$ are $10^{-14}$ away from 0, and all the entries of
  $UDV^\dagger-A$ are $3\cdot10^{-14}$ away from $0$, up to the accuracy of the
  multiplication, addition, and subtraction operations described in
  \ref{D:ieee754}. See the function \verb|VerifiedZSVD| in the accompanying code.
\end{dsc}
\begin{ntt}[Using the symmetry of the automorphism group]\label{N:group}
  One may use the symmetry group of the curve acting on the theta
  characteristics to save computations. Using the entire group is somewhat
  tricky from a computational point of view as the decomposition series has
  three terms, however, using just the
  semi-direct product of the bi-elliptic involutions and the rotations
  on the coordinates is reasonably easy: It is an index two subgroup of the
  automorphism group, and has a two term decomposition series.
  This is what we do in practice; see the functions \verb|group_mul| and
  \verb|group_act_on_theta| in the accompanying code. We denote this group by
  $G_0$.
\end{ntt}
\begin{ntt}
  For natural $J < 40$, let $q, j$ be the integer quotient and
  remainder of $J$ by 8, where $j$ is uniquely expressed as $j=j_0+2j_1+4j_2$
  for $j_0, j_1,j_2\in\{0, 1\}$. Denote by $pt_J$ the coordinate-wise rotation
  to the right by $q$ of
  \[
  (\ii (1-2j_0)a; (1-2j_1); \ii(1-2j_2); a; 0)\in\BP^4,
  \]
  (where $\ii, a$ are as in Proposition~\ref{P:wiman_odd}).
  Denote the theta characteristic comprised of $pt_a, pt_b, pt_c, pt_d$ by
  \verb|P[a, b, c, d]|. Once we enumerate theta characteristics in
  Proposition-Definition \ref{PD:enumerate_theta} below, we denote a set
  of theta characteristics whose numbers are $a,b,c,\ldots$ by \verb|[a,b,c,...]|.
\end{ntt}
\begin{prd}\label{PD:enumerate_theta}
  The following list is an enumeration of the theta characteristics of
  $O_{160}$:

  In the group $Aut(W)$ orbit of $q_{30}^\bullet$:
  {\tiny
\begin{verbatim}
0:P[8, 12, 37, 39]  1:P[9, 13, 33, 35]  2:P[10, 14, 32, 34]  3:P[11, 15, 36, 38]
4:P[5, 7, 16, 20]  5:P[1, 3, 17, 21]  6:P[0, 2, 18, 22]  7:P[4, 6, 19, 23]
8:P[13, 15, 24, 28]  9:P[9, 11, 25, 29]  10:P[8, 10, 26, 30]  11:P[12, 14, 27, 31]
12:P[21, 23, 32, 36]  13:P[17, 19, 33, 37]  14:P[16, 18, 34, 38]  15:P[20, 22, 35, 39]
\end{verbatim}
  }
  {\tiny
\begin{verbatim}
16:P[0, 4, 29, 31]  17:P[1, 5, 25, 27]  18:P[2, 6, 24, 26]  19:P[3, 7, 28, 30]
20:P[8, 12, 36, 38]  21:P[9, 13, 32, 34]  22:P[10, 14, 33, 35]  23:P[11, 15, 37, 39]
24:P[4, 6, 16, 20]  25:P[0, 2, 17, 21]  26:P[1, 3, 18, 22]  27:P[5, 7, 19, 23]
28:P[12, 14, 24, 28]  29:P[8, 10, 25, 29]  30:P[9, 11, 26, 30]  31:P[13, 15, 27, 31]
\end{verbatim}
  }
  {\tiny
\begin{verbatim}
32:P[20, 22, 32, 36]  33:P[16, 18, 33, 37]  34:P[17, 19, 34, 38]  35:P[21, 23, 35, 39]
36:P[0, 4, 28, 30]  37:P[1, 5, 24, 26]  38:P[2, 6, 25, 27]  39:P[3, 7, 29, 31]
40:P[8, 12, 33, 35]  41:P[9, 13, 37, 39]  42:P[10, 14, 36, 38]  43:P[11, 15, 32, 34]
44:P[1, 3, 16, 20]  45:P[5, 7, 17, 21]  46:P[4, 6, 18, 22]  47:P[0, 2, 19, 23]
\end{verbatim}
  }
  {\tiny
\begin{verbatim}
48:P[9, 11, 24, 28]  49:P[13, 15, 25, 29]  50:P[12, 14, 26, 30]  51:P[8, 10, 27, 31]
52:P[17, 19, 32, 36]  53:P[21, 23, 33, 37]  54:P[20, 22, 34, 38]  55:P[16, 18, 35, 39]
56:P[0, 4, 25, 27]  57:P[1, 5, 29, 31]  58:P[2, 6, 28, 30]  59:P[3, 7, 24, 26]
60:P[8, 12, 32, 34]  61:P[9, 13, 36, 38]  62:P[10, 14, 37, 39]  63:P[11, 15, 33, 35]
\end{verbatim}
  }
  {\tiny
\begin{verbatim}
64:P[0, 2, 16, 20]  65:P[4, 6, 17, 21]  66:P[5, 7, 18, 22]  67:P[1, 3, 19, 23]
68:P[8, 10, 24, 28]  69:P[12, 14, 25, 29]  70:P[13, 15, 26, 30]  71:P[9, 11, 27, 31]
72:P[16, 18, 32, 36]  73:P[20, 22, 33, 37]  74:P[21, 23, 34, 38]  75:P[17, 19, 35, 39]
76:P[0, 4, 24, 26]  77:P[1, 5, 28, 30]  78:P[2, 6, 29, 31]  79:P[3, 7, 25, 27]
\end{verbatim}
}
, and in the group $Aut(W)$ orbit of $q_{32}^\bullet$:
{\tiny
\begin{verbatim}
80:P[24, 25, 34, 37]  81:P[28, 29, 32, 39]  82:P[30, 31, 33, 38]  83:P[26, 27, 35, 36]
84:P[2, 5, 32, 33]  85:P[0, 7, 36, 37]  86:P[1, 6, 38, 39]  87:P[3, 4, 34, 35]
88:P[0, 1, 10, 13]  89:P[4, 5, 8, 15]  90:P[6, 7, 9, 14]  91:P[2, 3, 11, 12]
92:P[8, 9, 18, 21]  93:P[12, 13, 16, 23]  94:P[14, 15, 17, 22]  95:P[10, 11, 19, 20]
\end{verbatim}
  }
  {\tiny
\begin{verbatim}
96:P[16, 17, 26, 29]  97:P[20, 21, 24, 31]  98:P[22, 23, 25, 30]  99:P[18, 19, 27, 28]
100:P[24, 25, 35, 36]  101:P[28, 29, 33, 38]  102:P[30, 31, 32, 39]  103:P[26, 27, 34, 37]
104:P[3, 4, 32, 33]  105:P[1, 6, 36, 37]  106:P[0, 7, 38, 39]  107:P[2, 5, 34, 35]
108:P[0, 1, 11, 12]  109:P[4, 5, 9, 14]  110:P[6, 7, 8, 15]  111:P[2, 3, 10, 13]
\end{verbatim}
  }
  {\tiny
\begin{verbatim}
112:P[8, 9, 19, 20]  113:P[12, 13, 17, 22]  114:P[14, 15, 16, 23]  115:P[10, 11, 18, 21]
116:P[16, 17, 27, 28]  117:P[20, 21, 25, 30]  118:P[22, 23, 24, 31]  119:P[18, 19, 26, 29]
120:P[24, 25, 33, 38]  121:P[28, 29, 35, 36]  122:P[30, 31, 34, 37]  123:P[26, 27, 32, 39]
124:P[1, 6, 32, 33]  125:P[3, 4, 36, 37]  126:P[2, 5, 38, 39]  127:P[0, 7, 34, 35]
\end{verbatim}
  }
  {\tiny
\begin{verbatim}
128:P[0, 1, 9, 14]  129:P[4, 5, 11, 12]  130:P[6, 7, 10, 13]  131:P[2, 3, 8, 15]
132:P[8, 9, 17, 22]  133:P[12, 13, 19, 20]  134:P[14, 15, 18, 21]  135:P[10, 11, 16, 23]
136:P[16, 17, 25, 30]  137:P[20, 21, 27, 28]  138:P[22, 23, 26, 29]  139:P[18, 19, 24, 31]
140:P[24, 25, 32, 39]  141:P[28, 29, 34, 37]  142:P[30, 31, 35, 36]  143:P[26, 27, 33, 38]
\end{verbatim}
  }
  {\tiny
\begin{verbatim}
144:P[0, 7, 32, 33]  145:P[2, 5, 36, 37]  146:P[3, 4, 38, 39]  147:P[1, 6, 34, 35]
148:P[0, 1, 8, 15]  149:P[4, 5, 10, 13]  150:P[6, 7, 11, 12]  151:P[2, 3, 9, 14]
152:P[8, 9, 16, 23]  153:P[12, 13, 18, 21]  154:P[14, 15, 19, 20]  155:P[10, 11, 17, 22]
156:P[16, 17, 24, 31]  157:P[20, 21, 26, 29]  158:P[22, 23, 27, 28]  159:P[18, 19, 25, 30]
\end{verbatim}
  }
\end{prd}
\begin{proof}
  This is the content of the function \verb|fill_odd_theta| in the
  accompanying code. It applies the group action as in \ref{N:group} on
  the odd theta characteristics we got in Proposition~\ref{P:wiman_odd}.
\end{proof}
\begin{thm}\label{T:gaps}
  Then intersections of $O_{160}\choose 2$ with sets of pairs of Steiner system
  are distinct $G_0$ orbits which break as follows:
\ \\ \ \\
  \begin{tabular}{l c c c}
    number of orbits of this ``numerical'' type & 3 & 3 & 12\\
    size of orbit & 5 & 5 & 40\\
    number of pairs in the intersection ${O_{160}\choose 2}\cap \gS_\ga$&
    48&32&24.
  \end{tabular}
\end{thm}
\begin{proof}
  We first note that in order to see 18 orbit representatives one merely has
  to run the accompanying computer program, which also prints the data in the
  table above.
  The function \verb|build_all_steiner| in the program
  picks a representative $a:=\{\gth_1,\gth_2,\gth_3,\gth_4\}\in{O_{160}\choose 4}$
  for each $(\BZ/2)^4\rtimes \BZ/5$ orbit (as explained in \ref{N:group}),
  and performs on it the following computation:
  \begin{enumerate}
  \item\label{I:single2}
    First run the function \verb|maybe_two_k| in the accompanying code:
    Let $n\leq16$ be the number distinct points on $W\subset |K_\bW|^*$
    in the formal sum $\sum_{\gth\in a}\gth$.
    We compute the $n\times 12$ matrix $M$ from Corollary~\ref{C:single_svd},
    and it's SVD. The computation errors incurred are the one
    coming from evaluating representatives of the points $q_{ij}^k$
    (the error incurred here is in the floating point representation of
    $\sqrt{\phi}$), computing the degree 2 monomials of coordinates of said
    points which happens in the function \verb|fill_pts|, decomposing to
    irreducible representation as in
    Proposition~\ref{P:irrep_decomp} which happens in the function
    \verb|irrep_decomp|, and
    finally performing the SVD in Corollary~\ref{C:single_svd}.
    For all the $M$s we encounter, the computed
    singular values of $M$ lie in $[0,10^{-14}]\cup[10^{-2},10^3]$.
  \item\label{I:double2}
    Next, if $M$ has a computed singular value in the first segment above, and
    if the number of multiple points in $\sum_{\gth\in a}\gth$ is some $n_2>0$
    we run the function \verb|maybe_two_k_at_least_double| in the
    accompanying code:
    Let $M_2$ be the $2n_2\times k_1$ matrix from Corollary~\ref{C:svd_double},
    then we compute it's SVD. The computational errors we
    accumulate on top of
    the of the ones from the previous step come from the getting back from
    irreducible representations basis to monomial basis in the function
    \verb|anti_decomp|, computation of the
    gradient as in Remark~\ref{R:grad_comp} in the function
    \verb|q_grad_at_pt|,
    a dot product of two vectors of
    length $5$ per Lemma~\ref{L:tangent}, and the SVD from
    Corollary~\ref{C:svd_double} itself. For all the $M_2$s we encounter,
    the computed singular values of $M_2$ lie in $[0,10^{-13}]\cup[10^{-2},10]$.
  \item\label{I:triple2} Finally, if the singular values from \ref{I:double2} are
    in $[0, 10^{-13}]$, and $p$ is a triple point in $\sum_{\gth\in a}\gth$,
    then we run the function \verb|maybe_two_k_at_least_triple|:
    This function performs the verification in
    Proposition-Definition~\ref{PD:Hessian}, and in practice the computation
    in
    in \ref{D:half_grad_coef}. The computational errors accumulated on
    top of the ones we already have come from the operations in
    \ref{D:half_grad_coef}, and from a gradient computation as
    in Remark~\ref{R:grad_comp}. For all the computations we encounter,
    the difference corresponding to the one in \ref{triple} in Proposition
    \ref{P:whats2K}
    is in the set $[0, 10^{-14}]\cup[10^{-2}, 10]$.
  \end{enumerate}
  Let $\ol{A}$ be the set of $a$'s so that either
  $M$ has no singular value in the first segment in \ref{I:single2}, or it does
  but $a$ has double points, and $M_2$ has no singular value in the second
  segment in \ref{I:double2}, or it does, but the value computed in
  \ref{I:triple2} is in the second segment. Then by
  Proposition~\ref{P:whats2K},
  and by \ref{D:ieee754} and \ref{D:svd_accuracy}, for all $a\in\ol{A}$,
  $\sum_{\gth\in a}\gth\neq 2K_\bW$.
  The function \verb|build_all_steiner| then verifies that the conditions of
  Proposition~\ref{P:all2K} holds for $A:={O_{160}\choose 4}\sm \ol{A}$. Hence
  by Proposition~\ref{P:all2K}, the output of \verb|build_all_steiner| is
  a partition of $O_{160}\choose 2$ into it's intersection with the
  sets of pairs of Steiner systems of the curve $W$.
\end{proof}
\begin{prp}\label{P:Ic}
  For the Steiner systems in $6$ of the $12$ orbits of Steiner system described in the rightmost
  column of
  Theorem~\ref{T:gaps} (i.e.\ $240$ such systems),
  the space
  $\left(\spn\{q_{\{\gth,\gth+\ga\}}\}_{\gth\in\gS_{W,\ga}}\right)^\perp$ is 2 dimensional.
  The projections of this space on each of the 1st, 2nd, and 5th
  representations in Proposition~\ref{P:irrep_decomp} are non-trivial,
  and the projections on the 3rd and 4th representations is trivial.
\end{prp}
\begin{proof}
  In the function \verb|build_vc_alpha_Ic2_perp| in the accompanying code
  we show --- subject
  as usual to \ref{D:ieee754} and \ref{D:svd_accuracy} --- that for all partial
  Steiner systems, the singular values
  of the matrix whose rows are complex norm 1 representative of
  $\{q_{\{\gth,\gth+\ga\}}\}_{\gth\in\gS_{W,\ga}}$ are in $[0, 10^{-14}]\cup[10^{-2},10]$.
  Moreover, in the function \verb|handle_vc_a_plus_ic2_dim_13| we show that for
  all the $240$ partial systems above,
  $13$ of the values are in the second segment. Since the dimension
  cannot be bigger than $13$, it is exactly $13$.

  As for the projections on the irreducible representations, it is clear
  from Corollary~\ref{C:irrep_decomp} that the projections on the
  3rd and 4th representations are trivial. As for the other representations,
  we project, and compute (numerically, and subject again to \ref{D:ieee754})
  the magnitude is in the segment $[0.3, 1]$.
\end{proof}
\begin{cor}\label{C:Ic}
  The intersection of the $13$ dimensional $V_{W,\ga}$ is $I_2(W)$.
\end{cor}
\begin{proof}
  The intersection is orthogonal to the span of union of unitary orthogonal
  spaces to the $V_{W,\ga}$s, which, applying the symmetry group and the second
  part of Proposition~\ref{P:Ic}, is the
  direct sum of the 1st, 2nd, and 5th representations from
  Proposition~\ref{P:irrep_decomp}
\end{proof}
\begin{proof}[Proof of the main theorem]
  Follows from Theorem~\ref{T:gaps}, Corollary~\ref{C:Ic} and
  Proposition~\ref{P:specialize}
\end{proof}
Having proven the main theorem, we now move to some remarks about implementation
details:
\begin{dsc}[Testing the dedicated code]\label{D:tests}
  From a falsifiability point of view, a software test can merely prove that
  some code is wrong, which happens if and when the test fails. However,
  verifying full or intermediate results in ways which
  are easier to compute than, and independent of, the computations
  use to achieve the result give a heuristic affirmation to the claim that
  the code is performing its intended function. All the
  ``framework'' (i.e.\ OS, compiler, libraries) pieces behind the code used in
  the proof of \ref{T:gaps} and
  \ref{P:Ic} are extensively (to say the least)
  tested. There is one piece which is not: the dedicated code itself.
  The tests for this dedicated code are described below, and the claim therein
  is that they really are both independent from what they verify, and verify
  things which are far easier to verify than compute.
  \begin{itemize}
  \item The verification discussed in \ref{D:svd_accuracy}.
  \item In the function \verb|fill_pts|, we verify that the representatives of
    the points composing the $160$ theta
    characteristics evaluate to $0$ on the $Q_i$s, up to \ref{D:ieee754}.
  \item In the function \verb|fill_hyperplanes| we verify that the theta
    hyperplanes have double zeroes on all their intersection
    points with the curve, up to \ref{D:ieee754}.
  \item Moving from monomial lexicographic ordering to irreducible
    representations and back --- as in Proposition~\ref{P:irrep_decomp} --- are
    inverses of one another  up to \ref{D:ieee754}.
    (this is verified only in one direction, in the function \verb|anti_decomp|).
  \item The partial Steiner system we get satisfy
    Propositions~\ref{P:different_ij}, \ref{P:same_ij2} and
    Remark~\ref{R:gamma_orbits}:
    The numbers of sizes of partial Steiner system of each size are
    verified in the main
    program, the sizes of the conjugacy classes under the group action is
    verified in the function \verb|build_orbit_reps|,
    and the way they ``break'' between different $O_{ij}-O_{i',j'}$ is
    verified in the function \verb|verify_and_print_pair_structure|.
  \item In Theorem~\ref{T:gaps} we show that --- in Proposition \ref{P:all2K}
    notations --- the
    points {\em in} $A$ which do not have quadruple intersection with the
    curve satisfy the boundaries in
    Proposition~\ref{P:whats2K}; whereas in
    order to prove Proposition~\ref{P:all2K} we only have to test the points
    which are {\em not in} $A$. Moreover, the bound used for the points in $A$
    is order of magnitudes less than what is required by
    Proposition~\ref{P:whats2K}. This is verified in the respective functions
    described
    in the proof of the theorem.
  \item For any $R$ conjugacy class $a$ in A (in the sense of
    Proposition~\ref{P:all2K}),
    \[\dim\left(I_2(C)+\spn\{q_{\{\gth_1,\gth_2\}}\}_{\{\gth_1,\gth_2\}\in a}\right)\leq 13.\]
    This is verified in the main function.
  \end{itemize}
\end{dsc}
\begin{rmr}[Supplying a ``traditional'' witness for the computation]
  Verifying the decomposition of $O_{160}$ into Steiner systems is of
  course easier than finding it. However, we can do even better:
  In each (partial) Steiner system, we can find pairs of pairs of
  theta characteristics, such that the 16 points in $\BP^4$ involved are all
  distinct.
  If --- on a given Steiner system --- the graph of pairs of theta characteristics
  connected by this property has one connected component, a spanning
  tree of this graph is a witness for the collection of pairs of
  theta characteristics being a Steiner system, in the sense that each of the 16
  tuples of points need only verification in the sense of \ref{single} in
  Proposition~\ref{P:whats2K} with infinite accuracy (which can be done
  over finite fields and using the Chinese remainder theorem). Below
  is such a witness to one of the sets of pairs of the partial Steiner systems
  from Proposition~\ref{P:Ic}, which was found in the function
  \verb|print_single_tree|:

  {\tiny
\begin{verbatim}
{[0,9,22,70],[0,9,24,142],[0,9,25,82],[0,9,26,83],[0,9,27,143],[0,9,31,42],[0,9,88,134],
[0,9,94,111],[0,9,114,149],[0,9,130,154],[2,20,30,49],[2,24,49,142],[2,40,49,71],[2,49,91,92],
[2,49,108,132],[2,49,112,129],[2,49,150,152],[4,7,140,141],[4,20,30,140],[5,6,80,81],
[5,20,30,80],[8,20,30,62],[22,48,60,70]}.
\end{verbatim}
  }
  By Corollary~\ref{C:Ic} thi witness proves the theorem.
\end{rmr}

\end{document}